\RequirePackage{fix-cm}
\documentclass[smallextended]{svjour3}       
\smartqed  
\usepackage{graphicx}
\usepackage{amssymb,amsmath,amsfonts}
\usepackage[shortlabels]{enumitem}
\usepackage{float}
\usepackage{cite}
\newtheorem{algorithm}{Algorithm}

\newcommand*{\ee}{{\rm e}}
\newcommand*{\jj}{{\rm j}}

\newcommand{\paren}[1]{\left(#1\right)}

\newcommand{\norm}[1]{\left\| #1 \right\|}

\newcommand{\VV}[1]{\mathtt{Vec}\left(#1\right)}
\DeclareMathOperator{\Id}{Id}

\begin{document}

\title{Subpixel image reconstruction using nonuniform defocused images
	\thanks{This project has received funding from the ECSEL Joint Undertaking (JU) under grant agreement No. 826589. The JU receives support from the European Union's Horizon 2020 research and innovation programme and Netherlands, Belgium, Germany, France, Italy, Austria, Hungary, Romania, Sweden and Israel.}
}
\titlerunning{Subpixel image reconstruction using nonuniform defocused images}

\author{Nguyen Hieu Thao, Oleg Soloviev, Jacques Noom and Michel Verhaegen}


\institute{Nguyen Hieu Thao\at
	Delft Center for Systems and Control,
	Delft University of Technology, 2628CD Delft, The Netherlands. Department of Mathematics, School of Education, Can Tho University, Can Tho, Vietnam.
	\email{h.t.nguyen-3@tudelft.nl, nhthao@ctu.edu.vn}\\
	Oleg Soloviev\at
	Delft Center for Systems and Control,
	Delft University of Technology,
	2628CD Delft, The Netherlands.
	Flexible Optical B.V., Polakweg 10-11, 2288 GG Rijswijk, The Netherlands.
	\email{o.a.soloviev@tudelft.nl}\\
	Jacques Noom\at
	Delft Center for Systems and Control,
	Delft University of Technology,
	2628CD Delft, The Netherlands.
	\email{j.noom@tudelft.nl}\\
	Michel Verhaegen\at
	Delft Center for Systems and Control,
	Delft University of Technology,
	2628CD Delft, The Netherlands.
	\email{m.verhaegen@tudelft.nl}
}


\maketitle

\begin{abstract}
	
	This paper considers the problem of reconstructing an object with high-resolution using several low-resolution images, which are degraded due to \textit{nonuniform defocus effects} caused by angular misalignment of the subpixel motions.
	The new algorithm, indicated by the Superresolution And Nonuniform Defocus Removal (SANDR) algorithm, simultaneously performs the nonuniform defocus removal as well as the superresolution reconstruction.
	The SANDR algorithm combines non-sequentially the nonuniform defocus removal method recently developed by Thao \emph{et al.} and the least squares approach for subpixel image reconstruction.
	Hence, it inherits global convergence from its two component techniques and avoids the typical error amplification of multi-step optimization contributing to its robustness. Further, existing acceleration techniques for optimization have been proposed that assure fast convergence of the SANDR algorithm going from rate $\mathcal{O}(1/k)$ to $\mathcal{O}(1/k^2)$ compared to most existing superresolution (SR) techniques using the gradient descent method.
	An extensive simulation study evaluating the new SANDR algorithm has been conducted. As no algorithms are available to address the combined problem, in this simulation study we restrict the comparison of SANDR with other SR algorithms neglecting the defocus aberrations. Even for this case the advantages of the SANDR algorithm have been demonstrated.
	
	
\end{abstract}

\section{Introduction}\label{s:intro}

Image-based quality control is one of the important tools used during the manufacturing and end quality checks in semiconductor~\cite{HUANG20151}, automotive~\cite{Zhou2019}, and many other industries.
For Industry 4.0, requiring a fully automatized quality checks, the spatial resolution (size of the smallest feature that can be inspected) is the key factor that affects the overall efficiency and throughput of the control tool.
For high-quality imaging systems used in these tools, the spatial resolution is defined as the quotient between the pixel size and the magnification, and thus for a higher resolution, either a smaller pixel size or a larger magnification is required.
Larger magnification corresponds to smaller field of view (FOV in Fig.~\ref{fig:imagingscheme}), often it is desirable to have a smaller pixel size.
However, there are technological and design limits to the magnification and the smallest pixels that can be manufactured and/or used in these tools and thus a computational approach to increasing spatial resolution provides an interesting alternative. 
The reconstruction of an object with high-resolution from several low-resolution (LR) images capturing the object at subpixel-offset positions, called the \textit{SuperResolution} (SR) problem, has been studied for many decades \cite{PelKerSch87, UrGro92}.
A number of solution approaches have been proposed for the SR problem, including direct methods \cite{KimBosVal90, KimSu93} and iterative algorithms \cite{SauAll87, NguMilGol01, FarRobElaMil04, SroCriFlu07, TakKan08, LiHuGaoTaoNin10, ZhaYuaSheLi11, ZhaLiShiLin11, LagGhaHakRag16, WanLinDenAn17, HuaSunYanFanLin17}.

Superresolution reconstruction is possible if the LR images are registered for different subpixel-offset positions of the object.
In practice, shifting the object at subpixel scale can be a major challenge to the SR problem and gives rise to a number of important questions that need to be addressed.
Camera shake and motion blur induced by the shifts have been analyzed in \cite{BasBlaZis96, KanMil13}.
Inaccuracy of subpixel registration has been considered in \cite{LeeKan03, TakMilProEla09}.
Ideally, the shifting process should not cause any variations in the object orientation with respect to the camera. However, this is not always the case in practice and such deviations cause undesirable deterioration of the data images and thus the reconstruction.
Imprecise displacements with respect to the optical axis would introduce defocus blurs in the acquired images.
More challenging, the shifting process can induce rotational movements of the object causing \textit{nonuniform defocus effects} in the data.
To the best of our knowledge, the latter challenge has not been considered in the literature of SR by subpixel motions.

In this paper, we consider the problem of reconstructing a Superresolution Image using Nonuniform Defocused images, called the SIND problem.
Our consideration was primarily motivated by the inspection of wafers in semiconductor industry and the basic hypotheses are mainly inspired by its practical context, but the resulting solution is also scalable for similar applications of computer vision.
As an alternative to shifting the object, LR images can be registered using multiple cameras whose optical axes are typically at different directions towards the object.
This also results in \emph{nonuniform defocus blurs} in the acquired data, and the SIND problem covers this challenge as a special case with known and fixed blurs.

Solution approaches to the SIND problem should address three main tasks, including estimation of nonuniform defocus models, removal of nonuniform defocus effects, and reconstruction of an SR image.
Assuming that the nonuniform defocus models have been estimated, this paper is devoted to the last two tasks.
More specifically, we propose a new algorithm to simultaneously perform both Superresolution reconstruction And Nonuniform Defocus Removal (SANDR).
The SANDR algorithm combines the nonuniform defocus removal method recently developed in \cite{ThaOleJacMic21} and the least squares approach \cite{VerVer07} for subpixel image reconstruction but not in a sequential manner.
Hence, it inherits global convergence from its two component techniques and avoids the typical error amplification of multi-step optimization contributing to its robustness.
Further, existing acceleration techniques for optimization \cite{BecTeb09} 
have been proposed that assure fast convergence of the SANDR algorithm going 
from rate $\mathcal{O}(1/k)$ to $\mathcal{O}(1/k^2)$ compared to most existing SR techniques using the gradient descent method, where $k$ is the number of iterations.

As, to our knowledge, no algorithms are available to address the SIND problem, we demonstrate the advantages of the SANDR algorithm over other SR algorithms neglecting the defocus aberrations, see Sect. \ref{subs:comparison}.
It is important to mention that the Projected Gradient (PG) and the so-called Sequential Minimization (SM) algorithms reported along with SANDR in the numerical section are also considered for the SIND problem for the first time.
Hence, comparing the SANDR algorithm with them is not a goal of this paper.

\section{Problem formulation}\label{s:problem formulation}

\subsection{Superresolution by subpixel motions}\label{subs:SR}

\newcommand{\rblabel}[1]{\raisebox{-.2cm}{\makebox[0em][r]{\textbf{#1}\hspace*{0.5em}}}}
\begin{figure}[tb]
	\centering
	\includegraphics[width=0.95\linewidth]{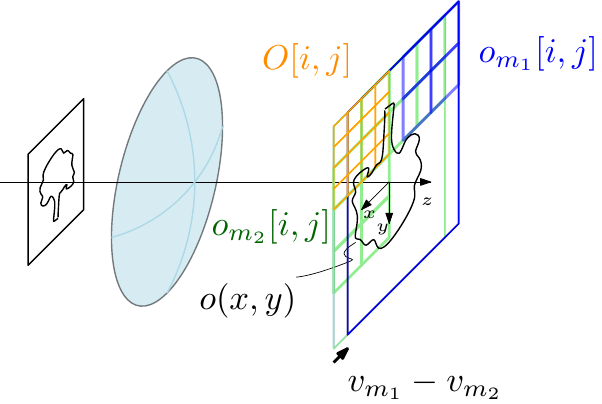}\vspace*{.25cm}\rblabel{a)}\\
	\includegraphics[width=0.95\linewidth]{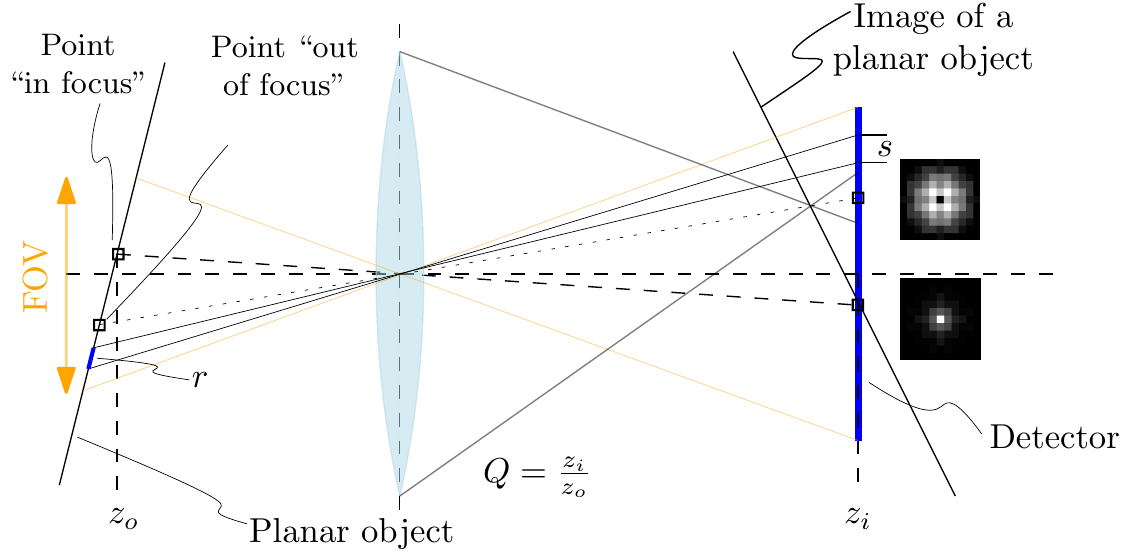}\vspace*{.25cm}\rblabel{b)}
	\caption{a) Imaging scheme used in the problem formulation.
		An imaging lens with magnification $Q$ creates image $o(x,y)$ of some planar object, which is registered by a detector with pixel size $s$ to obtain $M$ sampled images $o_m[i,j]$ (green, blue), with introduced subpixel offsets $v_m$ in each of them.
		The superresolution problem is to restore $O[i,j]$ representing $o(x,y)$ sampled with a finer grid (orange).
		b) Side view of the imaging scheme with an example of misalignment of the object and detector planes creating position-dependent defocus blur.}
	\label{fig:imagingscheme}
\end{figure}

Let $o_m$ be the LR images, created by an imaging system, that sample image $o$ of some planar object, see Fig.~\ref{fig:imagingscheme}.
Let each $o_m$ be registered with some subpixel offset $v_m$, with coordinates expressed in pixels of $o_m$.
We have the following sampling:
\begin{equation}\label{approx 1}
	o_m \sim T_{v_m}(o),\; (m=1,2,\ldots,M),
\end{equation}
where $T_v$ denotes the translation by a vector $v$.

Let $O$ be the SR image to be reconstructed. The ratio between the sizes of $O$ and $o_m$ is called the \textit{superresolution factor} and denoted by $\tau$.
In this paper, the images are assumed to be square and the superresolution factor is the same in both row and column directions for the sake of brevity.
As the shifts are measured in pixels of $o_m$, their coordinates with respect to $O$ should be scaled up by $\tau$.
Then we also have the following sampling:
\begin{equation}\label{approx 2}
	T_{\tau v_m}(O) \sim T_{v_m}(o),\; (m=1,2,\ldots,M).
\end{equation}

The combination of (\ref{approx 1}) and (\ref{approx 2}) leads to the following superresolution model: 
\begin{equation}\label{O&o_m}
	o_m \simeq D_{\tau}\circ T_{\tau v_m}(O),\; (m=1,2,\ldots,M),
\end{equation}
where $D_{\tau}$ is the \textit{downsampling operator} with rate $\tau$, see Sect. \ref{subs:downsampling operators} for the definition.

\begin{remark}[external blurs]\label{r:external blur}
	The imaging model (\ref{O&o_m}) can be extended as follows \cite{IraPel91, FarRobElaMil03, FarRobElaMil03.2, ParParKan03, SroCriFlu07, TakKan08, LiHuGaoTaoNin10, ZhaYuaSheLi11, ZhaLiShiLin11, LokSolSavVdo11, LagGhaHakRag16, WanLinDenAn17, HuaSunYanFanLin17}: 
	\begin{equation*}\label{with H}
		o_m \simeq D_{\tau}\circ H_m \circ T_{\tau v_m}(O),\; (m=1,2,\ldots,M),
	\end{equation*}
	where $H_m$ denote the external blurs often modelled as isoplanatic convolutions and assumed to be known.
	For our target application in wafer inspection, external blurs are not so relevant and thus left for brevity though they do not add major challenge to the problem under consideration.
	Instead, we handle the more challenging anisoplanatic blurs induced by angular misalignment of the shifts as detailed in the next section.
\end{remark}

\subsection{Nonuniform defocus effects}\label{subs:Anisoplanatic deconvolution}

In practice, the subpixel shifts can be accomplished by moving either the sample or the detector chip.
In both cases, some angular misalignment can be introduced, which can be difficult or costly (\textit{e.g.}, in terms of time/overall throughput) to eliminate completely.
Figure~\ref{fig:imagingscheme}b shows an example of a misaligned sample; a similar picture could be drawn for a misaligned detector, where geometrical distortions might also appear.\footnote{In this paper, we consider the geometrical distortions to be negligible compared to the position-dependent blur.} 
For simplicity, we do not discriminate between the object and detector misalignment.
Depending on the particular realisation and on the optical magnification of the system, this might create presence of position-dependent defocus blur in the image, which, as we show later, prevents the direct application of existing superresolution algorithms. 

We consider the challenge that the displacement process induces undesirable rotational movements of the object and the acquired images are degraded by nonuniform defocus blurs.
In this case, the theoretical LR images $o_m$ and the measured ones $i_m$ are related by
\begin{equation}\label{imaging model}
	i_m \simeq B_m(o_m),\; (m=1,2,\ldots,M),
\end{equation}
where the blur operators $B_m$ will be detailed shortly.

Let $o_m$ situate in $N$ defocus zones denoted by $D_n$ $(n=1,2,\ldots,N)$, for each of which the Point Spread Function (PSF) is modelled using the Fourier transform \cite{Goo05}:
\begin{equation}\label{p_n}
	p_n = \left|\mathcal{F}\paren{A \cdot \ee^{\jj d_n Z_2^0}}\right|^2\quad (n=1,2,\ldots,N),
\end{equation}
where the amplitude, product and square operations are elementwise, $\mathcal{F}$ is the two-dimensional Fourier transform, $A$ is the binary mask representing the camera aperture,\footnote{We assume that the diffraction-limited PSF is not resolved by the camera pixels.}
$\jj=\sqrt{-1}$ is the imaginary unit, $d_{n}$ is the (directional) distance from $D_n$ to the focal plane, and $Z_2^0$ is the Zernike polynomial of order two and azimuthal frequency zero.

We make use of the following model of nonuniform defocus blurs, whose physical relevance has been demonstrated, \textit{e.g.}, in \cite{ThaOleJacMic21}:
\begin{equation}\label{B_m_}
	B_m(o_m) = \sum_{n=1}^{N}\paren{\mu_{mn}\cdot o_m} * p_{n},\; (m=1,2,\ldots,M),
\end{equation}
where $*$ is the two-dimensional convolution, $\mu_{mn}$ are the mask functions of $o_m$ defined by: for $n=1,2,\ldots,N$,
\begin{equation*}\label{mu_mn}
	\mu_{mn}[i,j] =
	\begin{cases}
		1 & \text{if }\; o_m[i,j] \in D_n,
		\\
		0 & \text{otherwise}.
	\end{cases}
\end{equation*}

This paper considers planar objects and $d_{n}$ take the following form \cite[Sect. IIC,]{ThaOleJacMic21}:
\begin{equation}\label{d_n}
	d_{n} = d(n_0-n),\quad (n=1,2,\ldots,N),
\end{equation}
where $d$ is the \textit{Depth of Focus} (DoF) and $n_0$ is the \emph{focal position}.
The number of defocus zones $N$ defines the degree of defocus in an image and can be different for each $o_m$.
For brevity, it is taken the same in this paper.

The combination of (\ref{p_n}), (\ref{B_m_}) and (\ref{d_n}) yields the following blur model: ($m=1,2,\ldots,M$)
\begin{equation}\label{B_m}
	B_m(o_m) = \sum_{n=1}^{N}\paren{\mu_{mn}\cdot o_m} * \left|\mathcal{F}\paren{A\cdot \ee^{\jj d(n_0-n) Z_2^0}}\right|^2.
\end{equation}

\subsection{The SIND problem}\label{subs:SIND problem}

We consider the problem of reconstructing the Superresolution Image $O$ from the Nonuniform Defocused images $i_m$ according to the relations (\ref{O&o_m}) and (\ref{imaging model}), where $B_m$ are given by (\ref{B_m}).
It is referred to as the SIND problem.

\subsection{Optimization formulations}\label{subs:opt}

Combining (\ref{O&o_m}) and (\ref{imaging model}) yields the imaging model:
\begin{equation}\label{forward model}
	i_m = B_m\circ D_{\tau}\circ T_{\tau v_m}(O) + w_m,\quad (m=1,2,\ldots,M),
\end{equation}
where $B_m$ are given by (\ref{B_m}), and $w_m$ represent the discrepancies between the theoretical and the measured data, \textit{e.g.}, due to noise and model deviations.

In this paper, $w_m$ are assumed to be independent zero-mean random variables with jointly Gaussian-distributed entries.\footnote{This assumption is ubiquitous and it does not rule out the case of Poisson noise as the latter can be well approximated by a Gaussian distribution in view of the central limit theorem provided that the image is registered with a sufficiently large number of photon counts.}
For each $m=1,2,\ldots,M$, let $\mathcal{W}_m$ be the covariance matrix of $\VV{w_m}$, where $\mathtt{Vec}$ denotes the vectorization operator.
Then the \textit{maximum-likelihood} approach, \textit{e.g.}, \cite[Sect. 4.5.5,]{VerVer07}, applied to (\ref{forward model}) leads to the following minimization problem:
\begin{equation}\label{OP}
	\min_{O}\;\quad f(O) + \mathcal{G}(O),
\end{equation}
where $\mathcal{G}$ is the regularization capturing the 
physical attributes of $O$ (see Sect. \ref{subs:regularization}), and $f$ represents the data fidelity given by
\begin{align}\label{f(O)}
	f(O) &= \sum_{m=1}^{M}R_m(O)^T \mathcal{W}_m^{-1}\, R_m(O),
\end{align}
where $R_m$ ($m=1,2,\ldots,M$) are the fitting residual errors for the (blurred) LR images:
\begin{align*}
	R_m(O) &= \VV{B_m\circ D_{\tau}\circ T_{\tau v_m}(O)-i_m}.
\end{align*}

\begin{remark}[sequential optimization]\label{r:BIO}
	The residual error in (\ref{O&o_m}) is mainly due to the inaccuracy of the subpixel shifts while the one in (\ref{imaging model}) is more related to measurement noise and model deviations of (\ref{B_m}).
	When the latter is less severe than the former,\footnote{This is relevant to wafer inspection, where the camera is high-quality while inexactness of the subpixel shifts poses the major challenge.}
	one can also address (\ref{imaging model}) and (\ref{O&o_m}) sequentially via the following two-step optimization:
	\begin{equation}\label{BIO1}
		\min_{O}\;\sum_{m=1}^{M} S_m(O)^T {\Sigma}_m^{-1}\, S_m(O) + \mathcal{G}(O),
	\end{equation}
	where for $m=1,2,\ldots,M$,
	\begin{gather}\nonumber
		S_m(O) = \VV{D_{\tau}\circ T_{\tau v_m}(O)-\widehat{o}_m},
		\\\label{BIO2}
		\widehat{o}_m \in \arg\min_{o_m} Q_m(o_m)^T \mathcal{E}_m^{-1}\,Q_m(o_m) + \mathcal{H}\paren{o_m},
		\\\nonumber
		Q_m(o_m) = \VV{B_m(o_m)-i_m}.
	\end{gather}
	In the above, $\Sigma_m$ and $\mathcal{E}_m$ are respectively the covariance matrices representing the noise in (\ref{imaging model}) and (\ref{O&o_m}), and $\mathcal{H}$ is the regularization capturing the physical attributes of $o_m$.
	Sequentially minimizing (\ref{BIO2}) and (\ref{BIO1}) gives rise to the so-called SM algorithm (see Sect. \ref{subs:alg}), which suffers the typical error amplification of multi-step optimization compared to the proposed solution method for solving (\ref{OP}), see Sect. \ref{subs:solvability analysis}\&\ref{subs:noise analysis}.
\end{remark}

\subsection{Downsampling operators}\label{subs:downsampling operators}

The downsampling operator with integer rate $\tau$ is given by
\begin{equation}\label{D_tau}
	D_{\tau}([u]) = \frac{1}{\tau^2}\, \mathtt{conv}_{\tau} 
	\paren{u, \mathbf{1}_{\tau}},\quad (\forall\, u),
\end{equation}
where $\mathtt{conv}_{\tau}$ denotes the bivariate convolution operation with striking sizes $\tau\times \tau$,\footnote{The terminology is standard in the field of convolutional neural networks.} and $\mathbf{1}_{\tau}$ is the all-ones matrix of size $\tau\times \tau$.
The striking sizes define the size reduction in row and column directions.
$D_\tau$ produces only the average intensity value of every $\tau\times \tau$-block and hence it is not invertible without additional information of $u$.

\section{Solution approaches}\label{s:solution}

Solution approaches to the SIND problem should address three main tasks, including estimation of nonuniform defocus models, removal of nonuniform defocus effects, and reconstruction of an SR image.
Assuming that the nonuniform defocus models have been estimated, this paper is devoted to the last two tasks.
We first discuss regularization schemes for the SIND problem.

\subsection{Regularization functions}\label{subs:regularization}

SR methods often minimize a cost function consisting of data fidelity and regularization 
\cite{IraPel91, NguMilGol01, FarRobElaMil03, FarRobElaMil03.2, ParParKan03, 
	FarRobElaMil04, SroCriFlu07, TakKan08, LiHuGaoTaoNin10, ZhaYuaSheLi11, 
	ZhaLiShiLin11, LokSolSavVdo11, LagGhaHakRag16, WanLinDenAn17, 
	HuaSunYanFanLin17}.
Data fidelity is typically a norm of the residual between the theoretical and the measured data while regularization is driven by the \emph{a priori} known physical attributes of the solution.
The latter pertains to each particular application and is the main difference between existing SR techniques.
Total variation and Tikhonov regularization were considered in, \textit{e.g.}, \cite{NguMilGol01, FarRobElaMil03}.
The Bilateral Total Variation (BTV) was introduced in \cite{FarRobElaMil03} and later adapted in \cite{FarRobElaMil04, LiHuGaoTaoNin10, ZhaLiShiLin11, LagGhaHakRag16, WanLinDenAn17}.
In \cite{LagGhaHakRag16} BTV was used in combination with the Laplace operator while in \cite{LiHuGaoTaoNin10} it was used in combination with another regularization to enhance the consistence of the gradient variation.

In this paper, the images are assumed to have intensities in $[0,1]$, and the set of matrices satisfying this constraint is denoted by $\Omega$.
This constraint is easy to handle, but essential for the success of our proposed algorithms, where acceleration optimization mechanisms are exploited.
Its effectiveness has been widely known in the literature of deconvolution, see, \textit{e.g.}, \cite{WilSolPozVdoVer17, ThaOleJacMic21}. 
There are several approaches to this constraint, \textit{e.g.}, the penalty approaches using the associated distance function or its square.
In this paper, we make use of the \emph{indicator function} \cite{VA}:
\begin{equation}\label{indicator function}
	\iota_{\Omega}(x)
	= \begin{cases}
		0 & \text{ if } x\in \Omega,\\
		\infty & \text{ otherwise}.
	\end{cases}
\end{equation}

In our simulation results, this constraint is a precise regularization and hence its advantages over the other schemes are clearly observed, see Sect. \ref{subs:comparison}.

\subsection{The proposed algorithms}\label{subs:alg}

In view of Remark \ref{r:BIO}, the SIND problem can be addressed by solving (\ref{BIO2}) and (\ref{BIO1}) sequentially.
For each $m=1,2,\ldots,M$, (\ref{BIO2}) is the \textit{single-frame nonuniform defocus removal} problem recently studied in \cite{ThaOleJacMic21}.
Hence, it can be solved by the algorithm proposed in that paper, where its challenges including the typical ill-posedness were also discussed and global convergence of the proposed algorithm was also established. 
The main challenge of (\ref{BIO1}) is that the downsampling operator $D_{\tau}$ is not invertible, in particular, closed-form solutions for it are not available.
We propose to apply the regularization (\ref{indicator function}) and make use of the fast proximal gradient method introduced in \cite{BecTeb09}, often known as FISTA, for solving (\ref{BIO1}).

The algorithm resulted from this sequential approach will be referred to as the \textit{Sequential Minimization} (SM) algorithm for the SIND problem. However, we chose to skip its details for the sake of brevity.
The main advantages of the SM algorithm include its simplicity and the parallelism of (\ref{BIO2}) while its major disadvantage is the typical error amplification of multi-step optimization.

To overcome the drawback of SM, we next propose a new algorithm to simultaneously handle both Superresolution reconstruction And Nonuniform Defocus Removal (SANDR).
The SANDR algorithm combines the nonuniform defocus removal method developed in \cite{ThaOleJacMic21} and the least squares approach \cite{VerVer07} for subpixel image reconstruction but not in a sequential manner.
Hence, it inherits global convergence from its two component techniques and avoids the typical error amplification of multi-step optimization contributing to its robustness, see Sects. \ref{subs:solvability analysis}\&\ref{subs:noise analysis}.
Making use of the acceleration techniques for optimization of FISTA assures fast convergence of the SANDR algorithm going from rate $\mathcal{O}(1/k)$ to $\mathcal{O}(1/k^2)$ compared to most existing SR techniques using the gradient descent method.

For simplicity, the noise covariance matrices $\mathcal{W}_m$ in (\ref{f(O)}) are taken to be the identity matrix in the sequel. 
The repetitive term $(m=1,2,\ldots,M)$ following the subscript $m$ will be omitted for brevity.

In the sequel, $U_{\tau}$ will denote a right inverse of the downsampling operator $D_{\tau}$ defined in (\ref{D_tau}), \textit{i.e.}, $D_{\tau} \circ U_{\tau} =\Id$, the identity mapping.\footnote{$U_{\tau}$ is not unique and in general $D_{\tau}\circ U_{\tau}\neq U_{\tau}\circ D_{\tau}$.}
$U_{\tau}$ can be understood as a numerical upsampling operator, and in our numerical results, it is taken to be the interpolation with block constant values.
Recall that the translation by a vector with integer coordinates $(v_x,v_y)$, is given by
\begin{equation}\label{T}
	T_{(v_x,v_y)}(u)(r,c) = u(r-v_x,c-v_y),\quad (\forall\, u),
\end{equation}
where $(r,c)$ are the row-column coordinates of the pixels.

\begin{algorithm}[the SANDR algorithm]\label{al:SANDR}\ \\
	\emph{Input:} $i_m$ -- LR images, $B_m$ -- blur operators, $\lambda$ -- stepsize, $t^{(0)}$ -- initial acceleration stepsize, $K$ -- number of iterations, and $\varepsilon>0$.\\
	\emph{Initialization:} $X^{(0)} = O^{(0)}=\frac{1}{M}\sum_{m=1}^{M} T_{(-\tau v_m)}(U_{\tau}(i_m))$.\\
	\emph{Iteration process}: given $X^{(k)}$, $O^{(k)}$, $t^{(k)}$
	\begin{gather*}
		G_m^{(k)} = T_{(-\tau v_m)}\circ U_{\tau}\paren{\nabla f_m\paren{D_{\tau}\circ T_{\tau v_m}\paren{O^{(k)}}}},\\
		X_m^{(k+1)} = P_{\Omega}\paren{O^{(k)} - \lambda G_m^{(k)}},\\
		X^{(k+1)}=\frac{1}{M}\sum_{m=1}^{M}X_m^{(k+1)},\\
		t^{(k+1)} = \frac{1+\sqrt{1+4{t^{(k)}}^2}}{2},\\
		O^{(k+1)} = X^{(k+1)} + \frac{t^{(k)}-1}{t^{(k+1)}}\paren{X^{(k+1)}-X^{(k)}}.
	\end{gather*}
	\emph{Stopping criteria}: $k>K$ or
	\begin{equation}\label{stopping criterion}
		\sum_{m=1}^{M}\norm{G_m^{(k)}} >\sum_{m=1}^{M}\norm{G_m^{(k-1)}}+\varepsilon.
	\end{equation}
	\emph{Output:} $\widehat{O} = P_{\Omega}\paren{O^{(\mathtt{end})}}$.
\end{algorithm}

In Algorithm \ref{al:SANDR}, $P_{\Omega}$ is the projection operator associated with $\Omega$ and the functions $f_m$ are given by
\begin{equation*}
	f_m(x) = \frac{1}{2}\norm{B_m(x)-i_m}^2,\; (m=1,2,\ldots,M).
\end{equation*}

\section{Numerical simulations}\label{s:numerical}

As explained in Sect. \ref{subs:Anisoplanatic deconvolution}, a higher degree of defocus in an (LR) image corresponds to a larger number of defocus zones and smaller supports (nonzero entries) of the mask functions and vice versa.
To simplify simulation of random defocus levels in LR images, we chose to fix these parameters, but consider the DoF $d$ in (\ref{d_n}) as the single parameter quantifying the defocus in each image, called the \textit{blur coefficient} of the image in the sequel.
It is important to mention that our choice for convenience does not contradict the fact that DoF is a fixed physical parameter of the camera because underestimation of DoF does not introduce model deviations.\footnote{It only costs computational time as the number of defocus zones increases accordingly.}
The larger the blur coefficient is, the more the defocus blur in the image.

Simulation data is generated according to the forward imaging model 
(\ref{forward model}).
Except for the analysis regarding the number of input images in Sect. \ref{subs:number of shifts}, each data set consists of four images corresponding to the shift vectors $v_1=(0,0)$, $v_2=(1/2,0)$, $v_3=(0,1/2)$ and $v_4=(1/2,1/2)$.
Half of the images contain defocus blur varying in the vertical direction and half in the horizontal direction.
Unless otherwise specified, the common parameters are as in Table \ref{tbl:num_exp common params}.

\begin{table}[tb!]
	\caption{Parameters used in Sect. \ref{s:numerical}. $M$ is the number of LR images, $\mu$ -- size of the supports of mask functions in pixel rows/columns, $\rho$ -- PSF size (square), $\tau$ -- SR factor, $\lambda$ -- stepsize, and $t^{(0)}$ -- initial acceleration stepsize.}\label{tbl:num_exp common params}\vspace*{.25cm}
	\centering{
		\begin{tabular}[1\baselineskip]{r|cccccc}
			\textbf{Parameter} & $M$ & $\mu$ & $\rho$ & $\tau$ & $\lambda$ & $t^{(0)}$\\ \hline
			\textbf{Value} & 4 & 3 & 11 & 2 & 1 & 1\\ 
		\end{tabular}
	}
\end{table}

Except for the noise analysis in Sect. \ref{subs:noise analysis}, the data is corrupted with Poisson noise using the MATLAB function $\mathtt{imnoise}$.
The quality of SR reconstruction is measured by the \textit{Root Mean Square} (RMS) error of the restored SR image relative to the ideal one: ${\|\widehat{O}-O\|}\big{/}{\norm{O}}$.
The stopping criterion (\ref{stopping criterion}) is not implemented as it is not so relevant for simulations.

As no algorithms are available to address the SIND problem, we can only demonstrate its advantages over other SR methods neglecting the defocus effects.
It is important to mention that the Projected Gradient (PG) and the Sequential Minimization (SM) algorithms are also first considered for the SIND problem, and hence comparing the SANDR algorithm with them is not a goal of this section.
Instead, their own advantages and disadvantages in various problem settings will be of our primary interest.

\subsection{Comparison to known SR methods}\label{subs:comparison}

\begin{figure}[bt!]
	\centering
	\includegraphics[width=.95\linewidth]{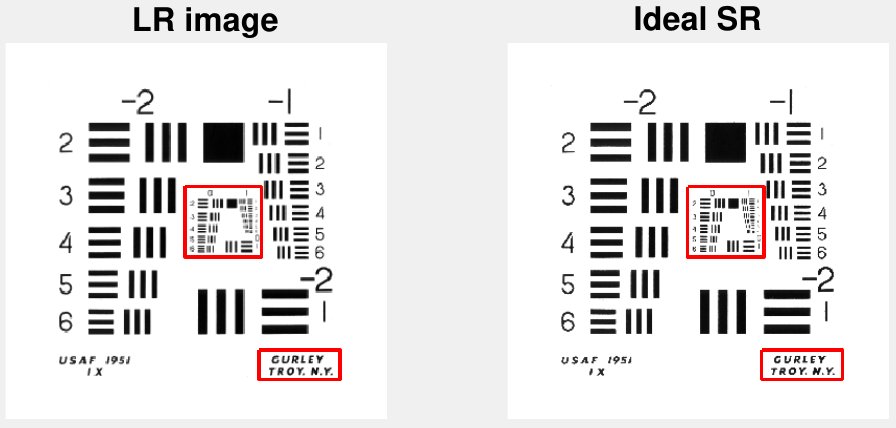}
	\caption{LR image (left) and the ideal SR (right). The ROIs are shown in Fig. \ref{fig:comparison_no_df ROIs} for visual comparison of different SR methods.}
	\label{fig:comparison_no_df data}
\end{figure}

\begin{figure}[bt!]
	\centering
	\includegraphics[width=.95\linewidth]{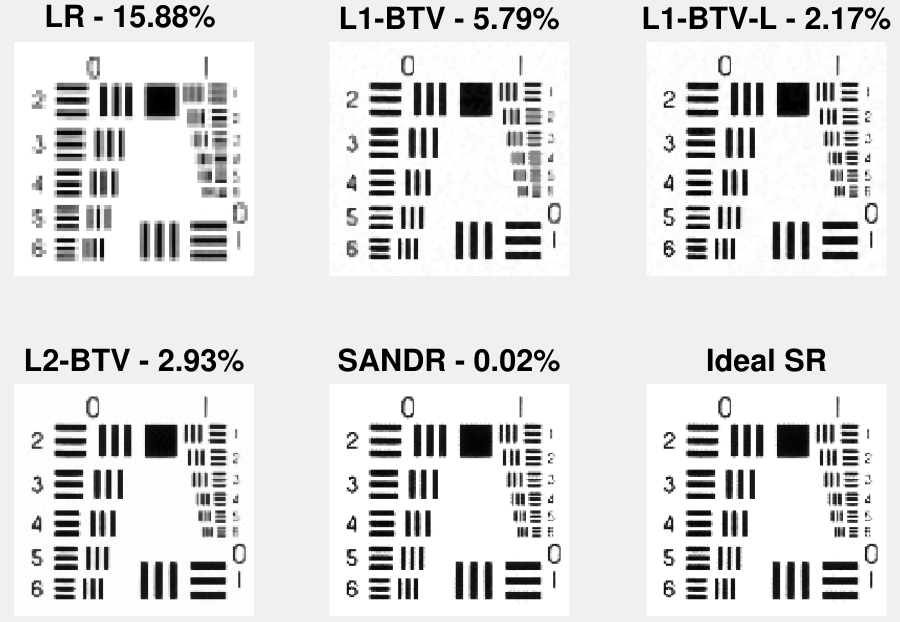}\vspace*{.25cm}\\
	\includegraphics[width=.95\linewidth]{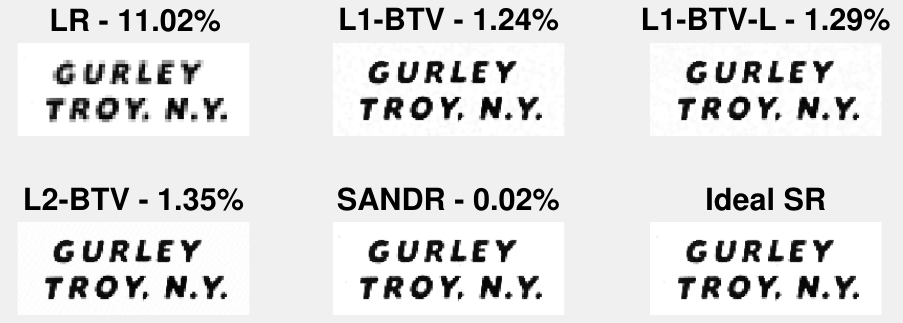}
	\caption{SR images obtained by L1-BTV, L1-BTV-L, L2-BTV and SANDR are shown together with an LR image and the ideal SR. Only the ROIs are shown for clarity. The relative RMS error of each ROI is also reported. The SANDR algorithm outperforms the other SR methods.}
	\label{fig:comparison_no_df ROIs}
\end{figure}

Most existing SR methods minimize a cost function consisting of data fidelity and regularization using the gradient descent method \cite{FarRobElaMil04, SroCriFlu07, TakKan08, LiHuGaoTaoNin10, ZhaYuaSheLi11, ZhaLiShiLin11, LagGhaHakRag16, WanLinDenAn17, HuaSunYanFanLin17}.
Data fidelity is typically the (weighted) $L_p$-norm ($1\le p\le 2$) of the residual between the theoretical and the measured data while regularization is driven by the \emph{a priori} known physical attributes of the solution, see Sect. \ref{subs:regularization}. In this section, we compare the SANDR algorithm with three existing SR methods minimizing (1) the $L_1$-norm with bilateral total variation (L1-BTV) \cite{FarRobElaMil04}, (2) the $L_1$-norm with BTV and Laplace operator (L1-BTV-L) \cite{LagGhaHakRag16}, and (3) the $L_2$-norm with BTV (L2-BTV) \cite{WanLinDenAn17}.
Each iteration of the algorithms is additionally followed by a projection on the constraint $\Omega$ to improve their performance, especially in terms of stability.
Note that without defocus effects, the SM and the SANDR algorithms coincide.

\begin{figure}[bt!]
	\centering
	\includegraphics[width=.95\linewidth]{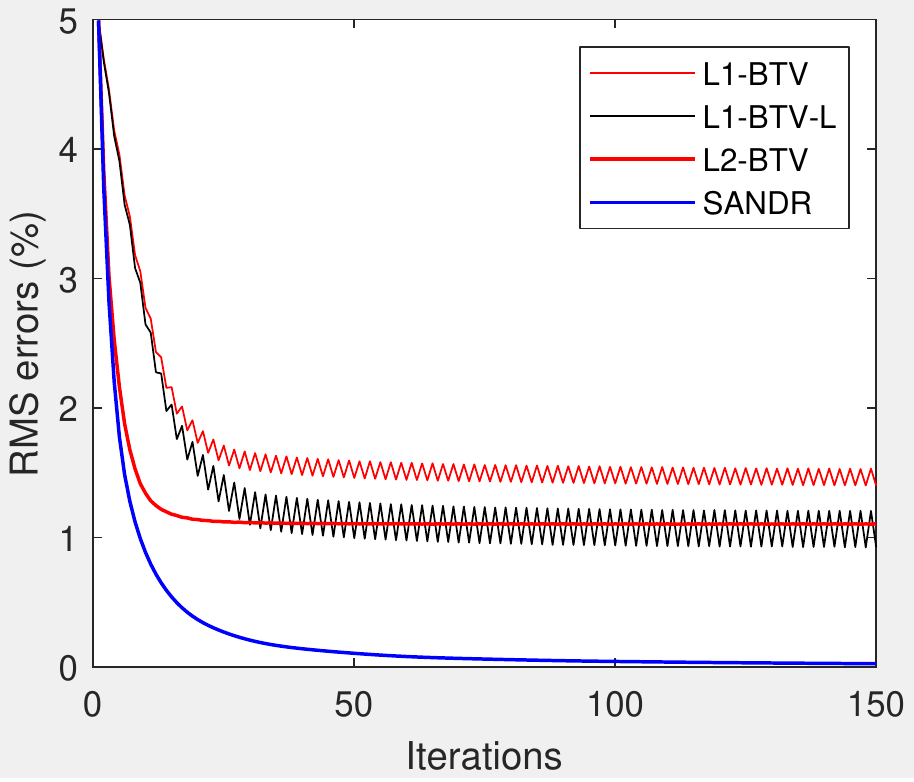}
	\caption{Relative RMS errors of the SR images obtained by L1-BTV, L1-BTV-L, L2-BTV and SANDR are shown in iterations. The SANDR algorithm is superior to the other methods in both convergence speed and accuracy.}
	\label{fig:comparison_no_df rrms}
\end{figure}

Figure \ref{fig:comparison_no_df data} shows an LR image (left) and the ideal SR (right).
The SR images obtained by L1-BTV, L1-BTV-L, L2-BTV and SANDR are shown in Fig. \ref{fig:comparison_no_df ROIs} together with an LR image and the ideal SR. Only the ROIs are shown for clarity.
The relative RMS error of each ROI is also reported.
The SANDR algorithm clearly outperforms the other methods both visually and in terms of RMS errors.
This is further explained in Fig. \ref{fig:comparison_no_df rrms}, where the RMS errors are shown in iterations.
The SANDR algorithm is far superior to the others in both convergence speed and accuracy.
Faster convergence is due to the acceleration feature of SANDR while higher accuracy can be explained by the fact that $\Omega$ is a precise regularization in this simulation problem.
Ripple behaviours of L1-BTV and L1-BTV-L in Fig. \ref{fig:comparison_no_df rrms} can be explained by the step-size being larger than the distance from the iteration to a local minimum.
This phenomenon is more likely to happen to $L_1$-norm cost functions since their gradient includes the sign function, which does not depend on the residual gap of the current iteration.\footnote{An advantage of $L_1$-norm cost functions is their potential to suppress outliers.}
Gradually decreasing the stepsize is a possible remedy for this issue, however, we chose not to distract the reader further in that direction because there does not exist a unified recipe for such tasks while the methods are not applicable to the SIND problem.

\subsection{Convergence properties}\label{subs:convergence properties}

\begin{figure}[bt!]
	\centering
	\includegraphics[width=.95\linewidth]{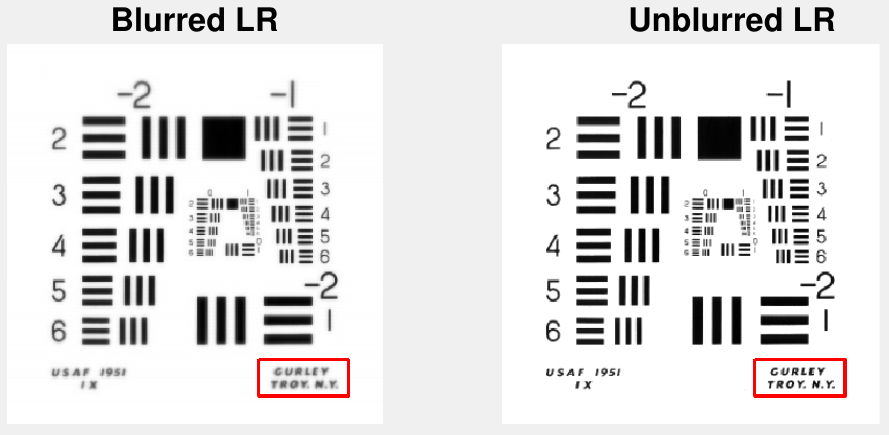}
	\caption{LR image with defocus varying in the vertical direction (left) and the unblurred one (right).}
	\label{fig:convergence data images}
\end{figure}

\begin{figure}[bt!]
	\centering
	\includegraphics[width=.95\linewidth]{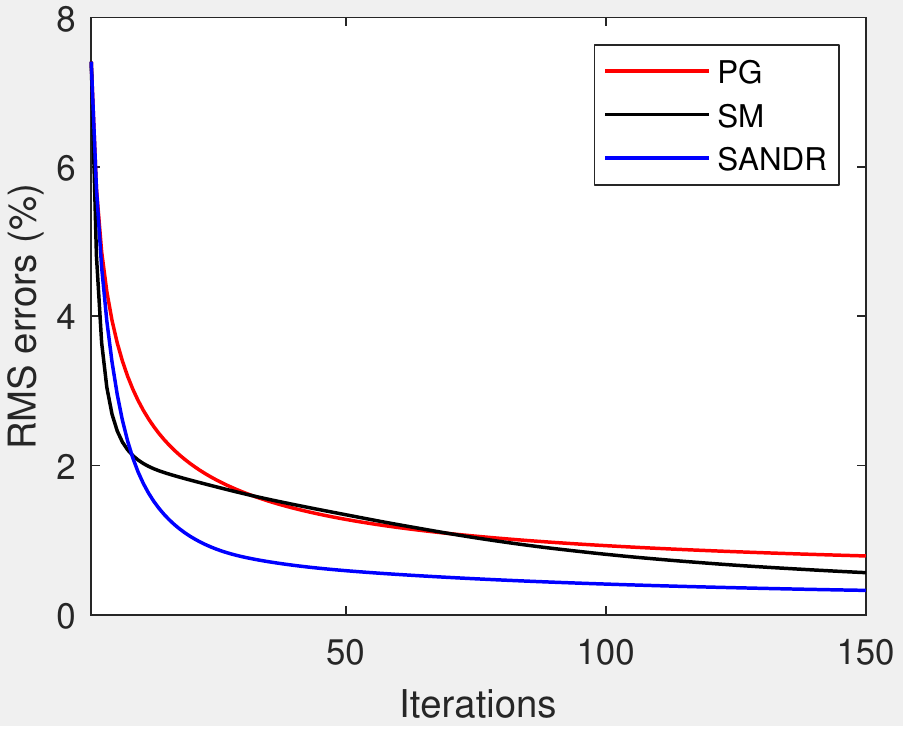}
	\caption{RMS errors of the SR images obtained by PG, SM and SANDR are shown in iterations. The algorithms exhibit convergence properties and without acceleration, PG (red) converges slower than SM (black) and SANDR (blue). The RMS error with 150 PG iterations is 0.79\%, corresponding to about 30 iterations of SANDR.}
	\label{fig:convergence rrms comparison}
\end{figure}

In this section, we demonstrate convergence properties of the SANDR algorithm along with the PG and SM methods.
We consider LR images of size $330\times 330$ pixels with $110$ defocus zones and blur coefficients randomly taken in the interval $[0.001, 0.06]$.
The other parameters are as in Table \ref{tbl:num_exp common params}.
One of the LR images with defocus effects varying in the vertical direction and its unblurred version are shown in Fig. \ref{fig:convergence data images}.

\begin{figure}[tb]
	\centering
	\includegraphics[width=.95\linewidth]{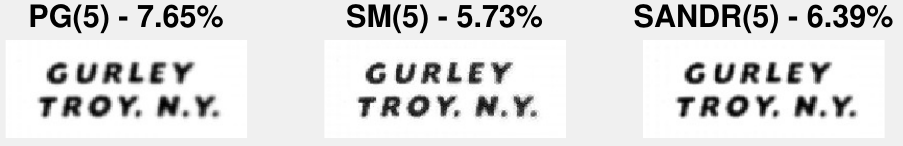}\vspace*{.25cm}\rblabel{a)}\\
	\includegraphics[width=.95\linewidth]{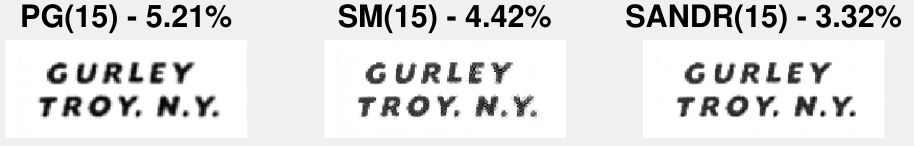}\vspace*{.25cm}\rblabel{b)}\\
	\includegraphics[width=.95\linewidth]{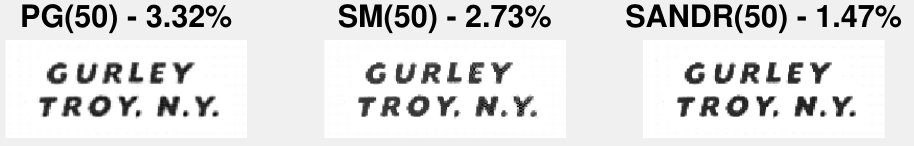}\vspace*{.25cm}\rblabel{c)}\\
	\includegraphics[width=.95\linewidth]{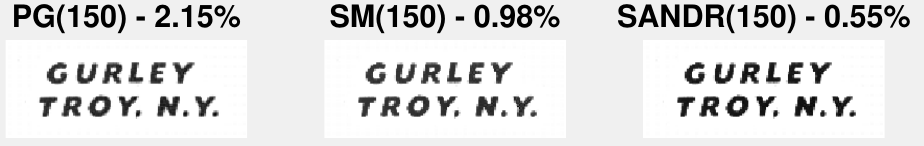}\vspace*{.25cm}\rblabel{d)}\\
	\includegraphics[width=.95\linewidth]{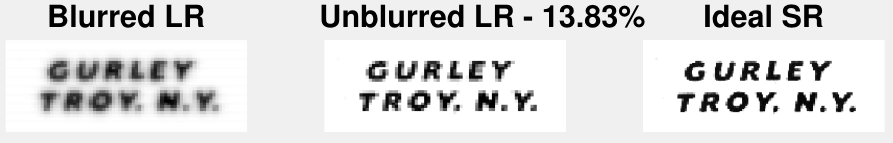}\vspace*{.25cm}\rblabel{e)}
	\caption{ROIs of the SR images obtained by a) 5, b) 15, c) 50, d) 150 
		iterations of PG, SM and SANDR are shown in comparison with e) the ones of 
		an LR image and the ideal SR. The RMS error of each ROI is also shown. The 
		RMS error of each unblurred LR image is around 13.83\%.}
	\label{fig:convergence ROI comparison}
\end{figure}

In Fig. \ref{fig:convergence rrms comparison} the RMS errors of the SR images obtained by PG, SM and SANDR are shown in iterations.
The algorithms exhibit convergence properties and without acceleration, PG (red) converges slower than SM (black) and SANDR (blue).
The RMS error with 150 iterations of PG is 0.79\%, corresponding to about 30 iterations of SANDR.
In Fig. \ref{fig:convergence ROI comparison} the ROIs of the SR images obtained by 5, 15, 50 and 150 iterations of PG, SM and SANDR are shown in comparison with the ones of an LR image and the ideal SR.
The relative RMS error of each ROI is also shown.
Note that the RMS error of each unblurred LR image is around 13.83\%.

\subsection{Solvability analysis}\label{subs:solvability analysis}

As explained at the beginning of Sect. \ref{s:numerical}, \textit{blur coefficients} quantify the defocus in the LR images.
The larger they are, the more challenging the SIND problem is.
In this section, we analyse the solvability of PG, SM and SANDR with respect to this parameter.
For each experiment, four blur coefficients of the LR images are randomly taken in the interval $[0.001, d_{\mathtt{max}}]$ with $d_{\mathtt{max}}$ ranging from 0.06 to 0.3.
The other parameters are as in Table \ref{tbl:df analysis_params}.\footnote{It is a trade-off between the computational complexity (number of iterations) and the restoration accuracy, in view of Figs. \ref{fig:comparison_no_df rrms} and \ref{fig:convergence rrms comparison} we chose to run $50$ iterations for each experiment.}

\begin{table}[bt!]
	\caption{Parameters used in Sect. \ref{subs:solvability analysis}. $d_{\mathtt{max}}$ is the upper bound of (random) blur coefficients, $N$ -- \# defocus zones, $n_0$ -- focal position, $(l,w)$ -- size of LR images, and $K$ -- \# iterations.}
	\label{tbl:df analysis_params}\vspace*{.25cm}
	\centering{
		\begin{tabular}[1\baselineskip]{r|ccccc}
			\textbf{Parameter} & $d_{\mathtt{max}}$ & $N$ & $n_0$ & $l=w$ & $K$\\
			\hline
			\textbf{Value} & 0.06 -- 0.3 & 55 & 28 & 165 & 50\\
		\end{tabular}
	}
\end{table}

\begin{figure}[bt!]
	\centering
	\includegraphics[width=.95\linewidth]{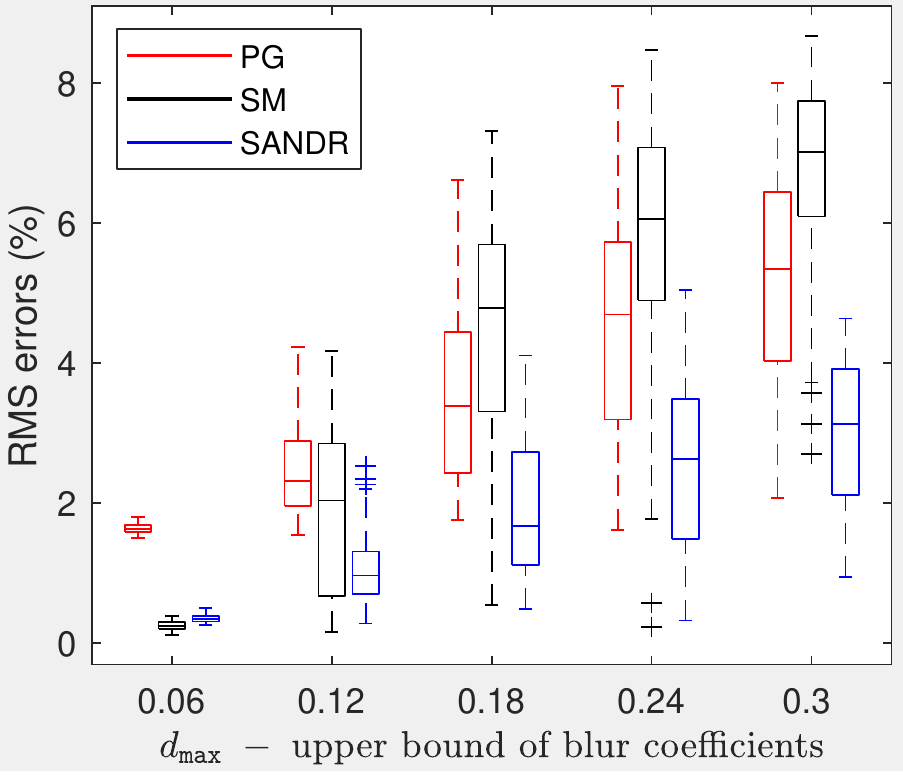}
	\caption{Solvability analysis of PG, SM and SANDR with respect to the distortion level of data images quantified by the \textit{blur coefficients}. One hundred experiments are reported for each value of $d_{\mathtt{max}}$ ranging from 0.06 to 0.3. The restoration errors increase for larger $d_{\mathtt{max}}$. SM works best for $d_{\mathtt{max}}$ up to 0.06, but it quickly becomes problematic for $d_{\mathtt{max}}$ from 0.12 due to its sequential optimization. SANDR is effective for blur coefficients up to 0.3 and outperforms PG and SM (for $d_{\mathtt{max}}\ge 0.12$) in both accuracy and consistency. The superiority becomes more significant for larger values of $d_{\mathtt{max}}$.}
	\label{fig:solvability analysis boxplot}
\end{figure}

\begin{figure}[bt!]
	\centering
	\includegraphics[width=.95\linewidth]{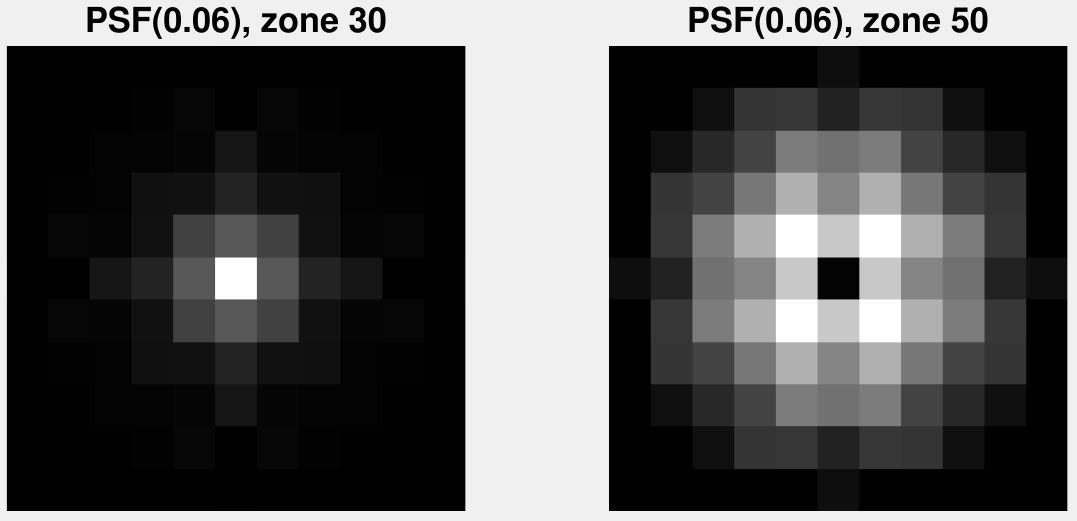}
	\caption{PSFs for the 30th (left) and 50th (right) defocus zones with blur coefficient 0.06.}
	\label{fig:solvability analysis PSFs}
\end{figure}

For each value of $d_{\mathtt{max}}$, one hundred experiments with PG, SM and SANDR are reported in Fig. \ref{fig:solvability analysis boxplot}, where the RMS errors of the obtained SR images with respect to the ideal one are presented.
The restoration errors increase for larger values of $d_{\mathtt{max}}$.
SM works best for $d_{\mathtt{max}}$ up to 0.06, but it quickly becomes problematic for $d_{\mathtt{max}}$ from 0.12 due to its sequential optimization.
SANDR is effective for blur coefficients up to 0.3.
It outperforms PG and SM (for $d_{\mathtt{max}}\ge 0.12$) in both accuracy and consistency, and the superiority becomes more significant for larger $d_{\mathtt{max}}$.
Higher accuracy is reflected by its smaller average restoration errors while more consistency is indicated by its smaller variances of the errors.
To visualize the \emph{blur coefficient} parameter, the PSFs for the 30th and 50th defocus zones (counted from the focal position) with blur coefficient 0.06 are shown in Fig. \ref{fig:solvability analysis PSFs}.
Recall that the distortion level of an image is proportional to the product of the blur coefficient and the zone position in view of (\ref{d_n}).

\subsection{Noise analysis}\label{subs:noise analysis}

\begin{figure}[bt!]
	\centering
	\includegraphics[width=.95\linewidth]{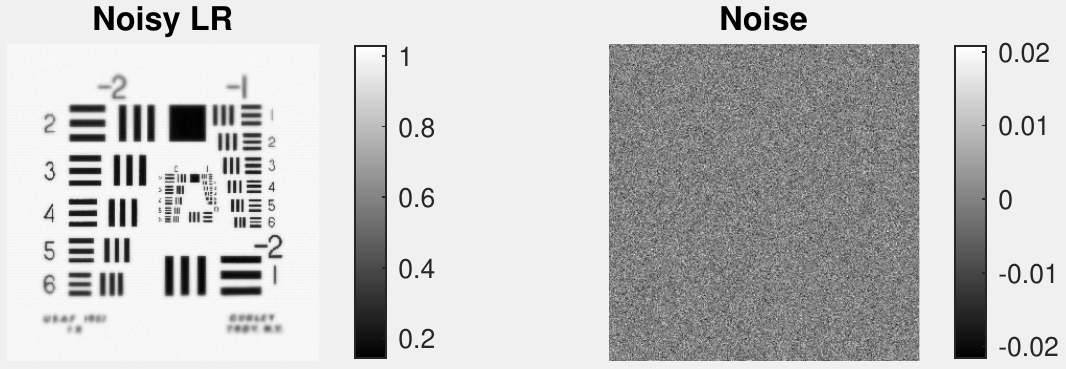}
	\caption{A noisy LR image (left) with SNR 45 dB and its residual (right) with respect to the noiseless one.}
	\label{fig:noisy image 45dB}
\end{figure}

We analyse the influence of noise on the performance of PG, SM and SANDR.
Five levels of Gaussian noise ranging from 45 to 65 dB (decibels) are considered.
Recall that the signal-to-noise ratio (SNR) in decibels is defined by: $\text{SNR} = 10\ln\left({P}/{P_0}\right)$, where $P$ and $P_0$ are the powers of the signal and noise, respectively.
To visualize the noise, an LR image with SNR 45 dB is shown in Fig. \ref{fig:noisy image 45dB} together with its residual relative to the noiseless one.

\begin{figure}[bt!]
	\centering
	\includegraphics[width=.95\linewidth]{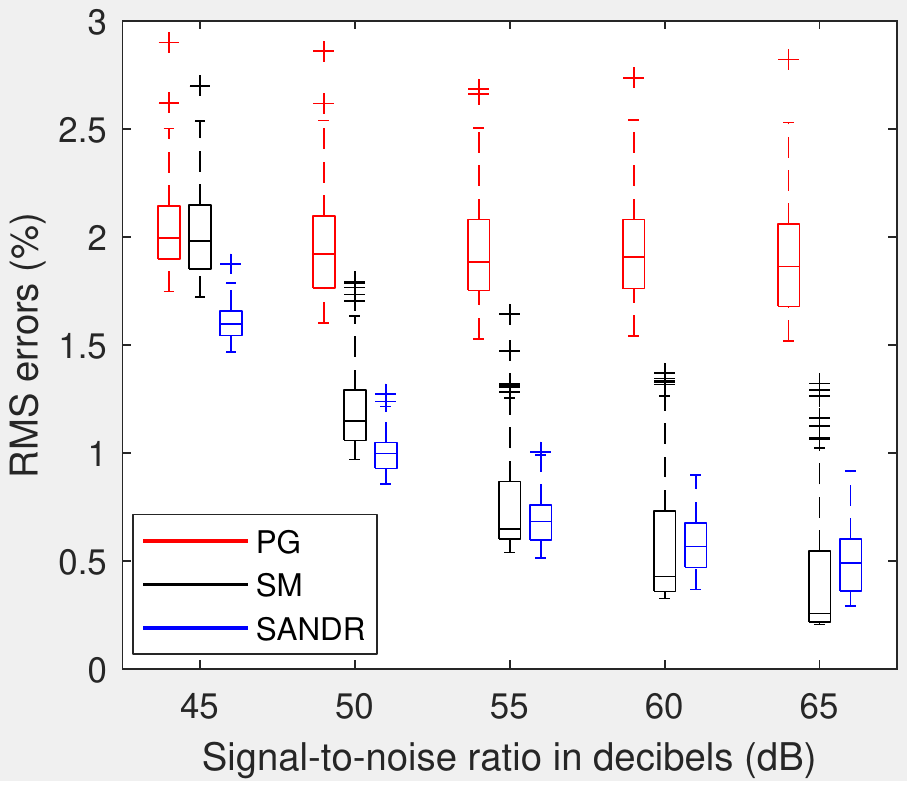}
	\caption{Noise analysis of PG, SM and SANDR. One hundred experiments are reported for each SNR from 45 to 65 dB. The reconstruction is more accurate for higher SNR. SM works best for SNR from 55 dB, but it quickly becomes problematic for SNR decreasing from 50 dB. Since SM and SANDR are acceleration variants of PG, they are less robust than PG.}
	\label{fig:noise analysis boxplot_df1u5}
\end{figure}

\begin{figure}[bt!]
	\centering
	\includegraphics[width=.95\linewidth]{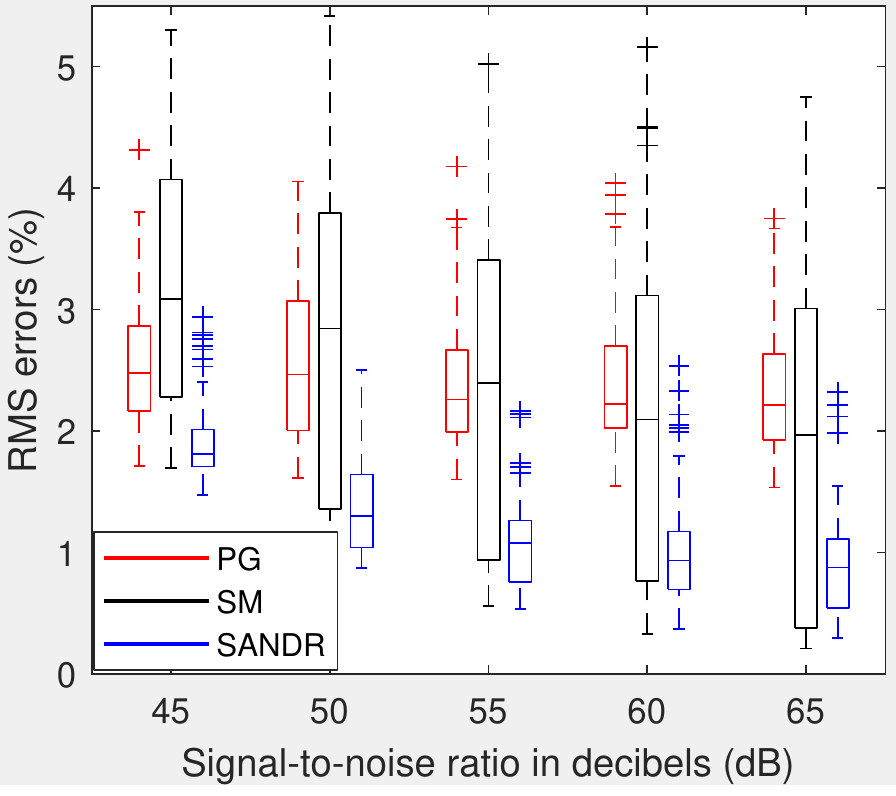}
	\caption{Experiments similar to those reported in Fig. \ref{fig:noise analysis boxplot_df1u5} but with $d_{\mathtt{max}}=0.12$ instead of $0.09$ show that SM deteriorates much faster than PG and SANDR for larger blur coefficients.}
	\label{fig:noise analysis boxplot_df2}
\end{figure}

For each SNR, one hundred experiments with $d_{\mathtt{max}}=0.09$ and the other parameters as in Table \ref{tbl:df analysis_params} are reported in Fig. \ref{fig:noise analysis boxplot_df1u5}, where the RMS errors of the SR images obtained by PG, SM and SANDR are presented.
The reconstruction is more accurate for higher SNR.
SM works best for SNR from 55 dB, but it quickly becomes problematic for SNR decreasing from 50 dB.
Its less robustness against noise compared to SANDR is due to its two-step optimization, see also Sect. \ref{subs:solvability analysis}.
It is not a surprise that SM and SANDR are less robust than PG because the former are acceleration variants of the latter and there is a typical trade-off between robustness and convergence speed.
In view of Fig. \ref{fig:noise analysis boxplot_df1u5}, it is worth thinking about PG for the SIND problem with SNR below 45 dB, but for higher SNR it is outperformed by the others.
It is important to recall that the conclusions drawn for the SM algorithm from Fig. \ref{fig:noise analysis boxplot_df1u5} are valid only for $d_{\mathtt{max}}$ up to 0.09, which seems to be a limit for it, see also Fig. \ref{fig:solvability analysis boxplot}.
To demonstrate this point, we do similar experiments but with slightly larger blur coefficients, $d_{\mathtt{max}}=0.12$ in place of $0.09$.
The results are summarized in Fig. \ref{fig:noise analysis boxplot_df2}, where SM deteriorates much more than PG and SANDR in comparison with Fig. \ref{fig:noise analysis boxplot_df1u5}.

\subsection{Number of input images}\label{subs:number of shifts}

The major practical challenge of SR by subpixel motions is to perform the shifts accurately.
Let us suppose that we are able to perform shifts at scale $1/\tau$ pixel, where $1<\tau\in \mathbb{Z}$.\footnote{In this study, inaccuracy of subpixel registration is subsumed in noise.}
Then there are at most $\tau^2$ LR images and one cannot expect to gain a SR factor greater than $\tau$.
In this section, we briefly study the influence of the number of input images on the quality of SR.
We consider $\tau=4$ and construct the SR image using 2, 4, 8 and 16 LR images, respectively.
In this experiment, $d_{\mathtt{max}}=0.06$ and the other parameters are as in Table \ref{tbl:df analysis_params}.

\begin{figure}[bt!]
	\centering
	\includegraphics[width=.95\linewidth]{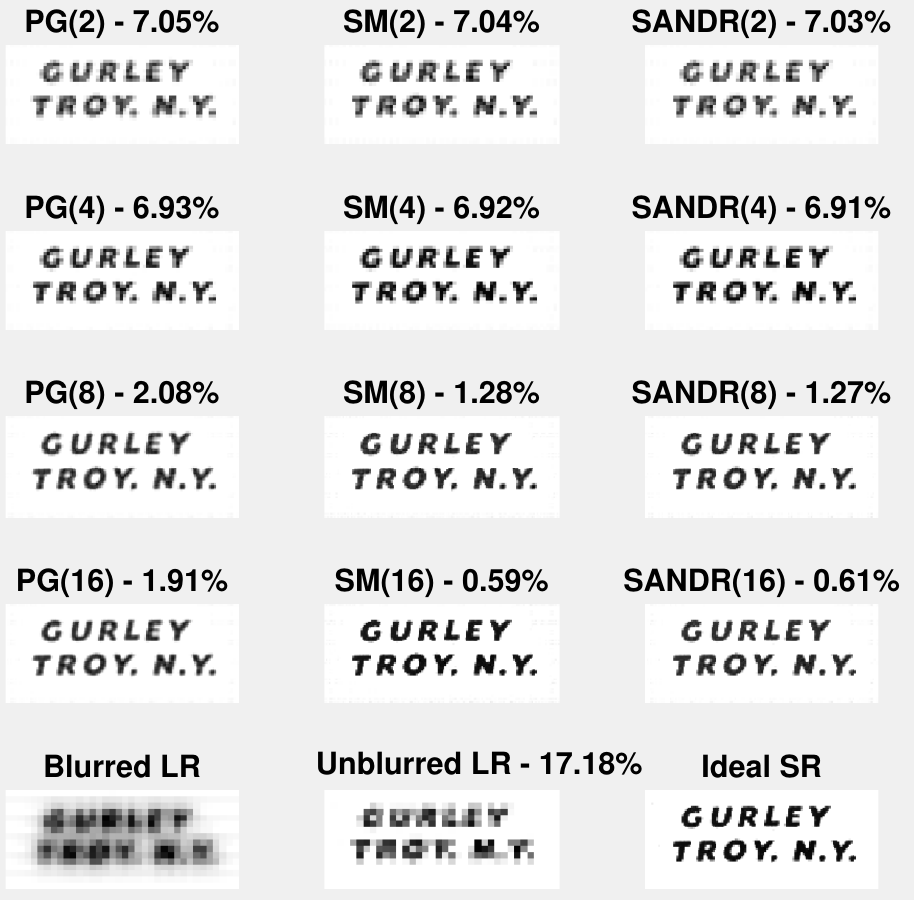}
	\caption{Influence of the number of input images on the SR images obtained by PG, SM and SANDR. ROIs of the SR images obtained with 2, 4, 8 and 16 LR images are shown together with their RMS errors. More input images result in higher quality of SR.}
	\label{fig:number of shifts}
\end{figure}

The numerical results are summarized in Fig. \ref{fig:number of shifts}, where only the ROIs and their RMS errors are shown for brevity.
It is clear that more input images result in higher quality of the SR and the observation is consistent for PG, SM and SANDR.
SR images obtained with two LR images (the first row) already shows improvement even in comparison with the unblurred LR images (the second in the last row).

\subsection{Image cropping}\label{subs:boundary effects}

Cropping the data images would introduce deviations to the imaging model (\ref{forward model}).
This issue does not arise in the previous sections since the simulation object there has almost constant intensity near the boundary. In this section, we study the influence of image cropping on the performance of the PG, SM and SANDR algorithms.

\begin{table}[tb!]
	\caption{Parameters used in Sect. \ref{subs:boundary effects}. $d_{\mathtt{max}}$ is the upper bound of blur coefficients, $N$ -- \# defocus zones, $n_0$ -- focal position, $(l,w)$ -- size of LR images, and $K$ -- \# iterations.}
	\label{tbl:boundary effects_params}\vspace*{.25cm}
	\centering{
		\begin{tabular}[1\baselineskip]{r|ccccc}
			\textbf{Parameter} & $d_{\mathtt{max}}$ & $N$ & $n_0$ & $l=w$ & $K$\\
			\hline
			\textbf{Value} & 0.09 & 85 & 43 & 255 & 50\\
		\end{tabular}
	}
\end{table}

\begin{figure}[bt!]
	\centering
	\includegraphics[width=.95\linewidth]{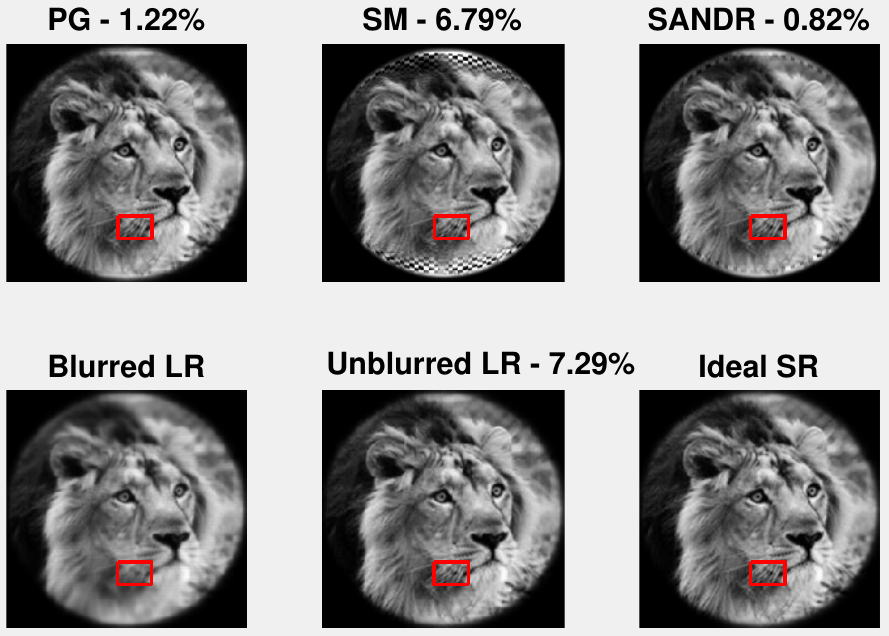}
	\caption{SR images obtained by PG, SM and SANDR are shown together with an LR image (bottom left) and the ideal SR (bottom right). The restoration error is smaller in the central region and becomes larger towards the boundary. The SANDR algorithm is the most effective for this problem while the SM algorithm is problematic due to high level of defocus effect. The SANDR and SM algorithms suffer more boundary effects than PG.}
	\label{fig:boundary effect}
\end{figure}

\begin{figure}[bt!]
	\centering
	\includegraphics[width=.95\linewidth]{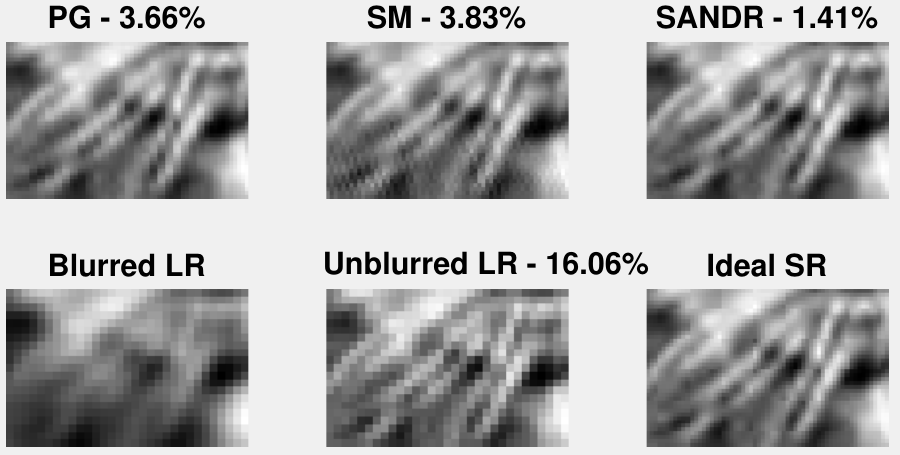}
	\caption{The ROIs marked in Fig. \ref{fig:boundary effect} are shown for a better inspection of finer details. Their RMS errors with respect to the ideal SR are also reported.}
	\label{fig:boundary effect ROI}
\end{figure}

Four images are generated according to (\ref{forward model}) with the parameters as in Table \ref{tbl:boundary effects_params}.
They are then windowed using Butterworth function to yield the LR images, one of which is shown at the bottom left of Fig. \ref{fig:boundary effect}.
The SR images obtained by PG, SM and SANDR are shown in Fig. \ref{fig:boundary effect} in comparison with an LR image and the ideal SR.
The SANDR algorithm is the most effective for this problem while the SM algorithm is problematic due to high level of defocus effect as discussed in Sect. \ref{subs:solvability analysis}.
The SM and SANDR algorithms suffer more boundary effects than PG since the former are more sensitive to noise than the latter as analyzed in Sect. \ref{subs:noise analysis}.
To reduce the restoration errors near the boundary, the Butterworth filter also need to be applied to every iteration update of the algorithms.
The RMS errors are computed for the central regions with 90\% in radius of the images.
The reconstruction error is smaller in the central region and increases towards the boundary.
The ROIs are zoomed out in Fig. \ref{fig:boundary effect ROI} for a better inspection of finer details.
The RMS error of each ROI is also reported.

\section{Concluding Remarks}\label{sec:remarks}

We have investigated the problem of constructing an object with high-resolution using several nonuniform defocused images, called the SIND problem.
Nonuniform defocus effects can arise in both standard techniques of data registration, including the use of multiple cameras and moving the object. However, the SIND problem has not been studied before.
We have proposed the efficient algorithm for SIND, called the SANDR algorithm, that can process both subpixel image reconstruction and nonuniform defocus removal simultaneously.
Important theoretical and practical aspects of the SANDR algorithm have been analyzed, including its global convergence, solvability, noise robustness, dependence on the number of LR images, and sensitivity to model deviations due to image croping.
We have demonstrated advantages of the SANDR algorithm over a number of existing superresolution methods without considering defocus effects because the latter cannot handle this additional challenge.
Our consideration was primarily motivated by the inspection of wafers in semiconductor industry, but the SANDR algorithm can be scalable for similar applications of computer vision in Industry 4.0.
\bigskip

\noindent\textbf{Funding.} This project has received funding from the ECSEL Joint Undertaking (JU) under grant agreement No. 826589. The JU receives support from the European Union's Horizon 2020 research and innovation programme and Netherlands, Belgium, Germany, France, Italy, Austria, Hungary, Romania, Sweden and Israel.
\bigskip

\bibliographystyle{plain}       
\bibliography{shortbib}   


\end{document}